\documentclass[12pt]{amsart}
\usepackage{amsthm}
\usepackage{amssymb}
\usepackage{amsfonts}
\usepackage{latexsym}
\usepackage{verbatim}
\usepackage{graphicx}
\usepackage{amsxtra}
\usepackage{enumerate}
\usepackage{comment}

\usepackage[usenames]{color}

\newtheorem{Corollary}{Corollary}
\newtheorem{Theorem}{Theorem}

\newtheorem{Lemma}{Lemma}[section]
\theoremstyle{definition}
\newtheorem{Definition}{Definition}[section]

\newtheorem*{Non-Planar Graph Theorem}{Non-Planar Graph Theorem}
\newtheorem*{Tucker's Theorem}{Tucker's Theorem}

\newtheorem*{Important Fact}{Important Fact}

\begin{document}

    \title{Asymmetric $2$-colorings of graphs}
    \date{\today}
    \author{Erica Flapan, Sarah Rundell, Madeline Wyse}
    


\address{Department of Mathematics, Pomona College, Claremont, CA 91711, USA}

\address{Department of Mathematics and Computer Science, Denison University, Granville, OH 43023, USA}

\address{Department of Near Eastern Studies, University of California, Berkeley, CA 94720, USA}

\thanks{The first author was supported in part by NSF Grant DMS-1607744}

\begin{abstract}  We show that the edges of every 3-connected planar graph except $K_4$ can be colored with two colors in such a way that the graph has no color preserving automorphisms.  Also, we characterize all graphs which have the property that their edges can be $2$-colored so that no matter how the graph is embedded in any orientable surface, there is no homeomorphism of the surface which induces a non-trivial color preserving automorphism of the graph.  
\end{abstract}

 \maketitle

\section{Introduction}

The study of asymmetric $2$-colorings of graphs embedded in $\mathbb{R}^3$ was originally motivated by the desire to classify the symmetries of non-rigid molecules.  While the symmetries of small molecules are induced by isometries of $\mathbb{R}^3$, large molecules have greater flexibility and hence some of their symmetries are not the result of a rigid motion. Such complex molecules can be represented by graphs in $\mathbb{R}^3$ where different colored edges represent different types of molecular chains or different types of bonds (see for example the representation of a molecular M\"{o}bius ladder in \cite{Si}).  Thus results about topological symmetries of colored graphs in $\mathbb{R}^3$ have potential applications to the study of symmetries of non-rigid molecules.

For example, Liang and Mislow \cite{LM1} used colored edges to distinguish between different molecular chains in their proof that certain families of proteins are chiral (i.e., topologically distinct from their mirror images).  In particular, after observing that these proteins all contained one particular embedding of the complete graph $K_5$ or the complete bipartite graph $K_{3,3}$ in $\mathbb{R}^3$, Liang and Mislow \cite {LM2} showed that by coloring some edges of these embedded graphs black and other edges grey, the colored graphs become topologically distinct from their mirror images.  They then conjectured that these colored graphs would remain topologically distinct from their mirror images even if they were embedded differently in $\mathbb{R}^3$.   

Motivated by this conjecture, Flapan and Li \cite{FL} proved that in fact the edges of any non-planar graph can be colored with two colors in such a way that every embedding of the graph in $\mathbb{R}^3$ is topologically distinct from its mirror image.  Furthermore they showed that, with the exception of the graphs $K_5$ and $K_{3,3}$, the edges of any non-planar $3$-connected graph can be $2$-colored so that for any embedding $\Gamma$ of the graph in $\mathbb{R}^3$,  no homeomorphism of $(\mathbb{R}^3,\Gamma)$ induces a non-trivial color preserving automorphism on $\Gamma$.

Our current results are stronger than those of Flapan and Li and apply to both planar and non-planar graphs.  To describe our results, we begin by introducing some terminology.  By a {\it graph} we shall mean a simple connected graph.  In particular, a graph does not have multiple edges joining the same pair of vertices nor loops containing only one vertex.  By a {\it $2$-coloring} of a graph $G$ we mean a coloring of each edge of $G$ with one of two colors.  Such a $2$-coloring is said to be \emph{intrinsically asymmetric} in a space $S$ if for any embedding $\Gamma$ of $G$ in $S$, every color preserving homeomorphism of the pair $(S,\Gamma)$ restricts to the trivial automorphism of $\Gamma$.  A $2$-coloring of $G$ is said to be \emph{intrinsically chiral} in an orientable space $S$ if for any embedding $\Gamma$ of $G$ in $S$, there is no orientation reversing color preserving homeomorphism of $(S, \Gamma)$.   We use the word \emph{intrinsic} in these contexts to emphasize that the asymmetry or chirality of the embedded graph depends only on the coloring of the graph $G$ and the space $S$ and not on the particular embedding of $G$ in $S$.  Using this terminology, Flapan and Li proved the following.
\medskip

\begin{Non-Planar Graph Theorem} \cite{FL}  \label{FlapLi}
{\it  \begin{enumerate}[(a)]
\item A non-planar 3-connected graph $G$ has a $2$-coloring that is intrinsically asymmetric in $\mathbb{R}^3$ if and only if $G$ is neither $K_{3,3}$ nor $K_5$.

\item Every non-planar graph $G$ has a $2$-coloring that is intrinsically chiral in $\mathbb{R}^3$. 
\end{enumerate}}
\end{Non-Planar Graph Theorem}

We now consider graphs which have the stronger property that their edges can be $2$-colored so they have no non-trivial color preserving automorphisms (independent of any embedding of the graph in $\mathbb{R}^3$).  We can describe this property more succinctly by building on terminology introduced by  Albertson and Collins \cite{AC}.  In particular, given a graph $G$, Albertson and Collins define the {\it distinguishing number} $D(G)$ to be the fewest number of colors needed to color the vertices of $G$ in such a way that $G$ has no non-trivial color preserving automorphisms.  We analogously define the  \emph{edge distinguishing number} $ED(G)$ of a graph $G$ to be the fewest number of colors needed to color the edges of $G$ so that $G$ has no non-trivial color preserving automorphisms.  Though their focus is on distinguishing numbers rather than edge distinguishing numbers, Albertson and Collins \cite{AC} make the observation that for all $n\geq 6$, $ED(K_n)=2$.

Observe that if a graph $G$ has $ED(G)\leq2$, then $G$ necessarily has a $2$-coloring which is intrinsically asymmetric in $\mathbb{R}^3$.  But one might wonder whether the converse is true.  In particular, can part (a) of the Non-Planar Graph Theorem be strengthened to show that every non-planar $3$-connected graph $G$ other than $K_5$ and $K_{3,3}$ actually has $ED(G)\leq 2$?  The following example shows that such a strengthening is not possible. 

 Let $H$ denote the graph consisting of $2^{7}+1$ copies of $K_5$ which are pairwise disjoint except along a triangle consisting of the edges $e_1$, $e_2$, and $e_3$ which all of the copies of $K_5$ share.   Observe that $H$ is non-planar and $3$-connected.  Now fix a 2-coloring of the edges $e_1$, $e_2$, and $e_3$, and note that for each $K_5$-subgraph of $H$, there are $2^7$ possible $2$-colorings of the seven edges of $K_5-\{e_1,e_2,e_3\}$.  Since $H$ contains $2^7+1$ subgraphs isomorphic to $K_5$, for any $2$-coloring of $H$ some pair of these subgraphs must have identical $2$-colorings.  Hence any $2$-coloring of $H$ will have a color preserving automorphism which interchanges a pair of identically colored $K_5$-subgraphs.  Recall however, that if we give $H$ the $2$-coloring from part (a) of the Non-Planar Graph Theorem, then for any embedding $\Gamma$ of $H$ in $\mathbb{R}^3$, no homeomorphism induces a non-trivial color preserving automorphism of $\Gamma$.  In particular, even though there is an automorphism which interchanges two identically colored $K_5$-subgraphs, no such automorphism can be induced by a homeomorphism of any embedding of $H$ in $\mathbb{R}^3$.

While we cannot strengthen part (a) of the Non-Planar Graph Theorem to show that every non-planar $3$-connected graph $G$ other than $K_5$ and $K_{3,3}$ has $ED(G)\leq 2$, in Section 2 we prove the following theorem which gives the desired result in the case of planar $3$-connected graphs.  

\begin{Theorem} \label{ED(G)} Let $G$ be a planar $3$-connected graph.  Then $ED(G) > 2$ if and only if $G=K_4$.  Furthermore, $K_4$ has no $2$-coloring which is intrinsically asymmetric in $\mathbb{R}^3$.
\end{Theorem}

Putting this theorem together with part (a) of the Non-Planar Graph Theorem, we obtain the following Corollary.

\begin{Corollary}  A $3$-connected graph $G$ has a $2$-coloring which is intrinsically asymmetric in $\mathbb{R}^3$ if and only if $G$ is not $K_4$, $K_5$, or $K_{3,3}$.
\end{Corollary}

In addition to considering symmetries of abstract graphs and graphs embedded in $\mathbb{R}^3$, we consider symmetries of graphs embedded in orientable surfaces.  In particular, for any graph $G$, we define the {\it surface edge distinguishing number} $SED(G)$ to be the smallest number of colors needed to color the edges of $G$ so that for any embedding of $G$ in any orientable surface $S$, no homeomorphism of $(S,G)$ induces a non-trivial color preserving automorphism of $G$.  Observe that for any graph $G$, we have $SED(G)\leq ED(G)$.

Tucker  \cite{Tu1,Tu} considered a concept that is related to the surface edge distinguishing number.  In particular, for a given surface $S$, Tucker classified all graphs $\Gamma$ embedded in $S$ such that for some $2$-coloring of $\Gamma$, no homeomorphism of $S$ induces a non-trivial color preserving automorphism on $\Gamma$.  Thus Tucker's result fixes an embedding of a given graph in a given surface and then chooses a $2$-coloring of the embedded graph that makes the embedding asymmetric.  By contrast, we are looking for a single $2$-coloring of an abstract graph $G$, which makes every embedding of the graph in every orientable surface asymmetric.

Before we state our result about surface distinguishing numbers, we introduce terminology to refer to some special types of graphs. The complete bipartite graph on partite sets of $n$ and $m$ vertices is denoted by $K_{n,m}$. The double star graph $S_{n,m}$ is  the graph obtained by connecting the vertex of degree $n$ in $K_{1,n}$ and the vertex of degree $m$ in $K_{1,m}$ by a path consisting of one or more edges.  Observe that $K_{1,1}$ is a single edge, $K_{1,2}$ is a path of 2 edges, and $S_{1,1}$ is a path of 3 or more edges.  Also, $C_n$ represents a cycle with $n$ vertices.  In Section 3, we prove the following theorem characterizing all graphs with surface edge distinguishing number greater than $2$.

\begin{Theorem}\label{SED(G)} A graph $G$ has $SED(G) > 2$ if and only if $G$ is a single edge, $C_3$, $C_4$, $C_5$, $K_4$, $K_5$, $K_{2,4}$, $K_{1,n}$ with $n\geq 3$, or $S_{n,m}$ with $n$ and $m$ odd and at least one of $n$ or $m$ greater than $1$.   \end{Theorem}

Finally, in Section 4 we characterize graphs which have a $2$-coloring which is intrinsically chiral in every orientable surface.  In particular, we prove the following.    

\begin{Theorem}\label{exceptions}  A graph $G$ has a $2$-coloring which is intrinsically chiral in every orientable surface if and only if $G$ has at least one vertex of degree at least $3$ and $SED(G)\leq 2$.    \end{Theorem}
\medskip

\section{Graphs with  $ED(G)> 2$}

In order to prove Theorem~\ref{K4Theorem}, we will use the following result of Tucker \cite{Tu1}.

\begin{Tucker's Theorem} \cite{Tu1}
\label{Tucker1}
Let $G$ be a graph and $A$ be a subgroup of $\mathrm{Aut}(G)$ such that no non-trivial element of $A$ fixes a pair of  adjacent vertices.  Then either the vertices of $G$ can be $2$-colored so that no non-trivial element of $A$ is color preserving or $G$ is one of the graphs $K_4$, $K_5$, $K_7$, $O_6$, or $O_8$.\end{Tucker's Theorem}

Note that the octahedral graph $O_n$ is defined as the graph obtained from $K_n$ by removing $\frac{n}{2}$ disjoint edges.  Thus all of the vertices of $O_n$ have degree $n-2$.  Observe that $O_6=K_{2,2,2}$ is the only graph in Tucker's Theorem which is planar and whose vertices have degree $4$.

\medskip

\begin{Lemma}\label{3edges}
Let $G$ be a graph embedded in an orientable surface $S$, and let $h$ be a homeomorphism of $(S, \Gamma)$ which fixes a vertex $v$ together with edges $e_1$, $e_2$, $e_3$ incident to $v$.  Then $h$ induces the identity automorphism on $G$.
\end{Lemma}

\begin{proof} Since $h$ fixes the vertices of $e_1$, $e_2$, $e_3$, $h$ cannot non-trivially rotate or reflect the edges incident to $v$.  Thus $h$ must fix every edge incident to $v$.  Also, there is a disk neighborhood $D$ of $v$ in $S$ which is setwise invariant under $h$. Define an orientation on $\partial D$ according to the order in which the edges $e_1$, $e_2$, and $e_3$ intersect $\partial D$.  This gives us an orientation on $S$ which is preserved by $h$.  

Now let $w$ be a vertex which is adjacent to $v$.  Then $h(w)=w$ and $h(\overline{vw})=\overline{vw}$.  Since $h$ preserves the orientation of $S$, $h$ cannot reflect the edges incident to $w$.  Also, since $h$ fixes $\overline{vw}$, $h$ cannot rotate the edges around $w$.  Thus $h$ fixes every edge incident to $w$.  Since $G$ is connected, we can inductively see that $h$ fixes every vertex of $\Gamma$.  
\end{proof}  

\medskip

We will use the following definition in the proof of Theorem~1.

\begin{Definition}
Let $G$ be a graph embedded in a surface $S$. We define the \emph{medial map} of $G$, denoted by $M(G)$, as the graph in $S$ obtained by placing a vertex, $v_e$ in the interior of each edge $e$ of $G$, and placing an edge between the vertices $v_{e_1}$ and $v_{e_2}$ if the edges $e_1$ and $e_2$ are adjacent on the boundary of a face in $S-G$.
\end{Definition}

\setcounter{Theorem}{0}
\begin{Theorem}\label{K4Theorem}
Let $G$ be a planar $3$-connected graph.  Then $ED(G) > 2$ if and only if $G=K_4$.  Furthermore, $K_4$ has no $2$-coloring which is intrinsically asymmetric in $\mathbb{R}^3$.\end{Theorem}

\begin{figure}[htb]
\centerline{\includegraphics[width=1.5in]{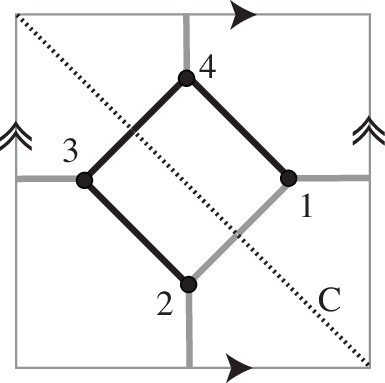}}
\caption{$(12)(34)$ is induced by a reflection of the torus through $C$.}
\label{K4T2}
\end{figure}

\begin{proof} Suppose that $G\not =K_4$.  By Whitney's Theorem \cite{Wi}, since $G$ is planar and $3$-connected, $G$
has a unique embedding in $S^2$.  Thus we abuse notation and refer to this embedded
graph also as $G$. We begin with some observations about the medial map $M(G)$.  By its construction, since $G$ is embedded in $S^2$, $M(G)$ is also embedded in $S^2$; and by the uniqueness of the embedding of $G$ in $S^2$ we know that $M(G)$ also does not depend on a particular embedding of $G$ in $S^2$.  Furthermore, since $G$ is 3-connected, $M(G)$ cannot have multiple edges with the same pair of vertices or loops through a single vertex, and hence $M(G)$ is a graph according to our definition.  Also, each edge of $G$ is contained in
precisely two faces of $S^2-G$.  Finally, each face of $S^2-G$ containing a given edge $e$ has exactly two
edges which are adjacent to $e$.  Thus as illustrated in Figure~\ref{4valent}, every vertex of $M(G)$ has degree $4$. 

\begin{figure}[htb]
\centerline{\includegraphics[width=.70\textwidth]{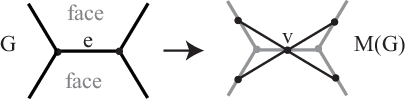}}
\caption{Every vertex of $M(G)$ has degree 4.}
\label{4valent}
\end{figure}

 Let $A=\mathrm{Aut}(G)$.  We know from Whitney's Theorem \cite{Wi} that every element of $A$ is induced by a homeomorphism of $(S^2,G)$.  Furthermore, since every homeomorphism of $(S^2,G)$ induces an automorphism of $M(G)$, we abuse notation and also consider $A$ as a subgroup of $\mathrm{Aut}(M(G))$.

  Now consider an automorphism $\alpha\in A$ which fixes a pair of adjacent vertices $x$ and $y$ of $M(G)$.  Let $a$, $b$, and $c$ be the vertices of $G$ such that $x$ is in the interior of $\overline{ab}$ and $y$ is in the interior of $\overline{ac}$ as illustrated in Figure~\ref{adjacent}.  It now follows that $\alpha$ fixes $a$, $b$, and $c$.   
   
\begin{figure}[htb]
\centerline{\includegraphics[width=.25\textwidth]{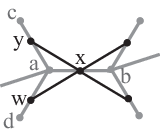}}
\caption{Every vertex of $M(G)$ has degree 4.}
\label{adjacent}
\end{figure} 

 Since the edge $\overline{ab}$ is in the boundary of precisely two faces of $S^2-G$, there is a unique vertex $w$ of $M(G)$ which is adjacent to $x$ and is in the interior of some edge $\overline{ad}$ of $G$ with $d\not =c$ (see Figure~\ref{adjacent}).   Now since $x$, $y$, $a$, $b$, and $c$ are fixed by $\alpha$, the vertex $w$ must also be fixed by $\alpha$.  It follows that $\alpha$ fixes the edges $\overline{ab}$, $\overline{ac}$, and $\overline{ad}$ of $G$.  It now follows from Lemma~\ref{3edges} that $\alpha$ is the identity automorphism of $G$.  Hence $\alpha$ is also the identity automorphism of $M(G)$.  Thus the identity is the only element of $A$ that fixes a pair of adjacent vertices of $M(G)$.

We can now apply Tucker's Theorem to conclude that since $M(G)$ is planar and all of its vertices have degree $4$, either $M(G)=K_{2,2,2}$ or
there is a $2$-coloring of the vertices of $M(G)$ such that no non-trivial automorphism in $A$ is color preserving. If $M(G)=K_{2,2,2}$, then $G$ would be $K_4$, which is contrary to our assumption.  Thus there is a $2$-coloring of the vertices of $M(G)$ such that no non-trivial automorphism in $A$ is color preserving.  This $2$-coloring of the vertices of $M(G)$ gives us a $2$-coloring of the edges of $G$ which has the property that no non-trivial automorphism in $\mathrm{Aut}(G)$ is color preserving.  Thus  $ED(G)\leq 2$.

\medskip

To show that $ED(K_4)>2$, we only need to consider $2$-colorings of $K_4$ with up to three grey edges, since a $2$-coloring with more than three grey edges is equivalent to the coloring obtained by interchanging grey and black edges.  All non-trivial $2$-colorings of $K_4$ are displayed in Figure~\ref{K4Chart} along with a non-trivial color preserving automorphism for each.  It thus follows that $ED(K_4)>2$. 
 
\begin{figure}[h!]
\centerline{\includegraphics[width=\textwidth]{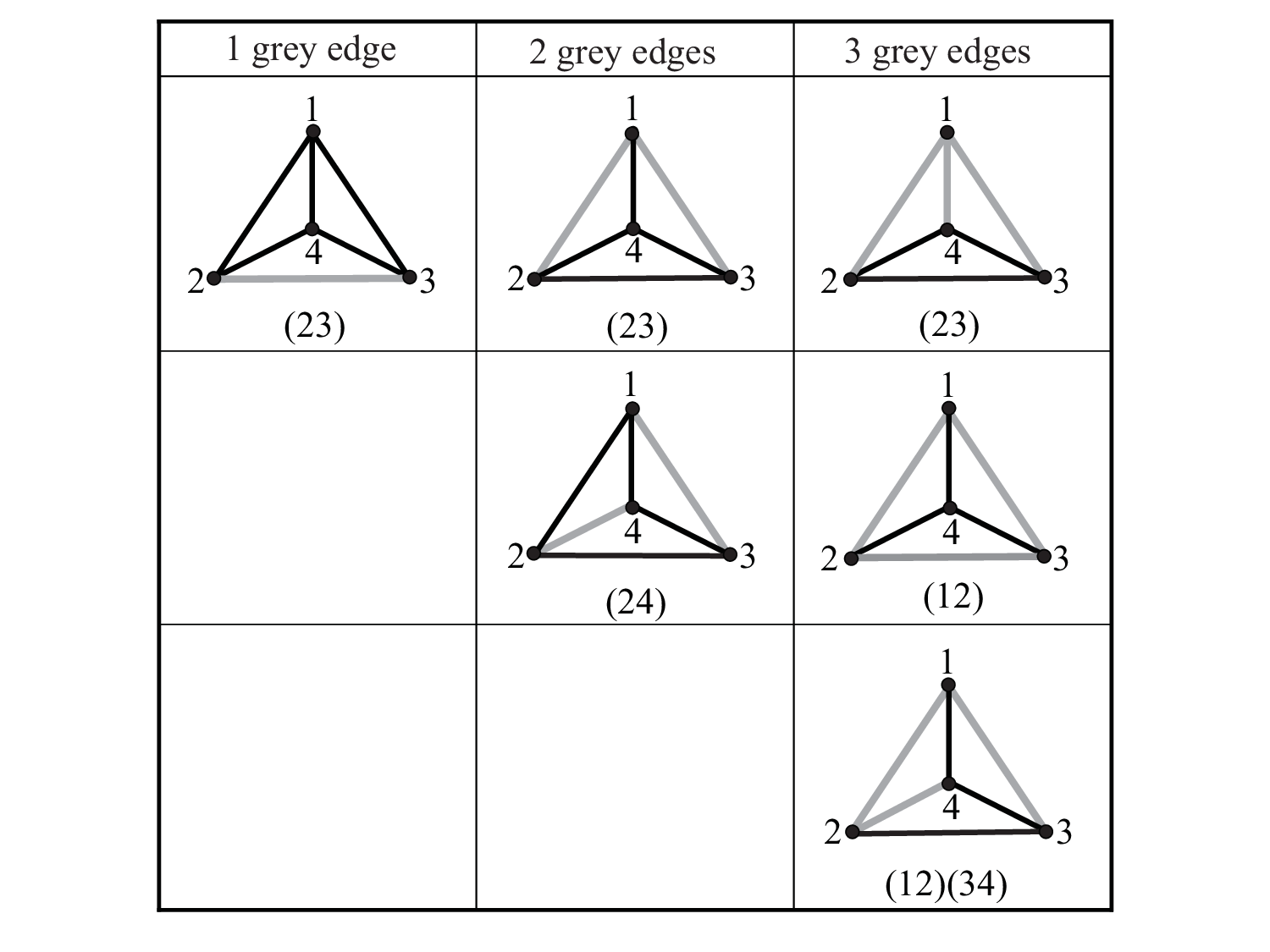}}
\caption{Every $2$-coloring of $K_4$ has a non-trivial color preserving automorphism.}
\label{K4Chart}
\end{figure}

In order to prove that no $2$-coloring of $K_4$ is intrinsically asymmetric in $\mathbb{R}^3$, we observe that the $2$-colorings of $K_4$ illustrated in Figure \ref{K4Chart} have no crossings, and hence can be viewed as embeddings of $K_4$ in $S^2$.  Now each of the automorphisms in the first two rows of Figure \ref{K4Chart} is induced by a color preserving reflection of $S^2$, which extends to a reflection of $\mathbb{R}^3$.  

For the $2$-coloring in the third row of Figure \ref{K4Chart}, we let $\Gamma$ denote the embedding of $K_4$ in a torus $T^2$  illustrated in Figure~\ref{K4T2}.  Then the automorphism $(12)(34)$ is induced by a reflection $f$ of $T^2$ through the $(1,1)$-curve $C$.  Now we embed $T^2$ as an unknotted torus in the $3$-sphere $S^3$ such that $f$ is induced on $T^2$ by a rotation $h$ of $S^3$ around $C$.  Let $p$ be a point on $C$ which is disjoint from the embedded graph $\Gamma$.  Now $h$ restricts to a homeomorphism of $\mathbb{R}^3=S^3-\{p\}$ which induces $(12)(34)$ on $\Gamma$.  \end{proof}

\medskip

\section{Graphs with $SED(G) > 2$}

 We now consider the surface edge distinguishing number.  We begin by proving three lemmas which will be used in the proof of Theorem~\ref{SED(G)}.
 
\begin{Lemma}
\label{RefereeLemma1}
Let $G$ be a graph containing a triangle as a proper subgraph.  If $SED(G) > 2$, then $G$ is $K_4$ or $K_5$.
\end{Lemma}

\begin{proof}  Let $v$, $x$, and $y$ be the vertices of a triangle in $G$.  Then without loss of generality, we can assume that $G$ has a vertex $z$ adjacent to $v$ which is distinct from $x$ and $y$.  We color the edges of the path $\overline{xyvz}$ black and color all of the other edges of $G$ grey.  Then any color preserving automorphism of $G$ must either fix each vertex in $\overline{xyvz}$ or interchange $x$ and $z$ as well as $y$ and $v$.  If an automorphism fixes each vertex in $\overline{xyvz}$, then it fixes each of the edges $\overline{xv}$, $\overline{yv}$, and $\overline{vz}$.

Now since $SED(G)>2$, there is an embedding $\Gamma_1$ of $G$ in some orientable surface $S_1$ such that a homeomorphism $h_1$ of $(S_1,\Gamma_1)$ induces a non-trivial color preserving automorphism of $\Gamma_1$.  By Lemma \ref{3edges}, we know that $h_1$ must interchange vertices $x$ and $z$ as well as $y$ and $v$.  Since $v$ is adjacent to $x$, this means that $y$ is adjacent to $z$ 

Next we consider a different $2$-coloring of $G$.  This time we color the path $\overline{xvzy}$ black and color the rest of the edges of $G$ grey.  Using an argument analogous to the above, we see that there is an embedding $\Gamma_2$ of $G$ in an orientable surface $S_2$ such that some homeomorphism $h_2$ of $(S_2,\Gamma_2)$ interchanges vertices $x$ and $y$ as well as $z$ and $v$.  Now since $v$ is adjacent to $y$, this means that $x$ is adjacent to $z$.  Thus $G$ contains the $K_4$-subgraph with vertices $x$, $y$, $z$, and $v$.

Now suppose that $G$ contains an additional vertex $p$.  Since $G$ is connected, without loss of generality $p$ is adjacent to $v$.  We now repeat the above argument with $p$ in place of $z$ to show that $p$ is adjacent to all of the neighbors of $v$.  Thus $G$ contains the $K_5$-subgraph with vertices $x$, $y$, $z$, $v$, and $p$.

By arguing inductively, we see that if $G$ has $n$ vertices then $G$ must be $K_n$.  However, we know from Albertson and Collins \cite{AC} that for all $n\geq 6$, $ED(K_n)=2$.  Thus  $G$ must be either $K_4$ or $K_5$.  \end{proof}\medskip

\begin{Lemma}
\label{RefereeLemma2}
Let $G$ be a graph which does not contain a triangle and which has a vertex $a$ of degree at least $3$ adjacent to at least two vertices each having degree at least $2$.  If $SED(G) > 2$, then $G=K_{2,4}$. 
\end{Lemma}

\begin{proof} Let $1$, $2$, and $3$ denote vertices that are adjacent to vertex $a$; and suppose that vertices $1$ and $2$ have degree at least $2$.  Then vertex $1$ is adjacent to some vertex $c$ which is distinct from $a$.  We color the path $\overline{c1a2}$ and the edge $\overline{a3}$ black and color the rest of $G$ grey.  Since $SED(G) > 2$, we know from Lemma~\ref{3edges} that $G$ must have an automorphism interchanging vertices $2$ and $3$.  By repeating this argument with vertex $1$ in place of $2$, we see that $G$ has an automorphism interchanging vertices $1$ and $3$.  Similarly, every pair of vertices adjacent to $a$ can be interchanged by some automorphism of $G$.  It follows that all of the neighbors of $a$ have the same degree.

Next suppose for the sake of contradiction that some neighbor of $a$ is not adjacent to $c$.  Without loss of generality assume that vertex $2$ is not adjacent to $c$.  Thus vertex $2$ is adjacent to a vertex $y$, which is distinct from $a$ and $c$.  Note that since $G$ does not contain a triangle, neither $c$ nor $y$ can be adjacent to $a$.  Now color $e=\overline{1c}$ and $P=\overline{3a2y}$ black and color the rest of the edges of $G$ grey, and let $f$ be a color preserving automorphism of $G$.  Then $f(e)=e$ and $f(P)=P$.  Furthermore, since $1$ is adjacent to $a$, and $c$ is not adjacent to either $a$ or $2$, we must have $f(1)=1$ and $f(c)=c$.  Now since $1$ is adjacent to $a$ but not to $2$, we must have $f(a)=a$.  From this it follows that $f$ fixes every vertex on $P$, and hence fixes three edges incident to $a$.   But now Lemma \ref{3edges} implies that this $2$-coloring of $G$ is intrinsically asymmetric in any orientable surface in which $G$ embeds.  As this contradicts our hypothesis that $SED(G) > 2$, we conclude that every neighbor of $a$ is adjacent to $c$.  Furthermore, by interchanging the roles of $a$ and $c$, we see that every neighbor of $c$ is adjacent to $a$.

Let $n$ denote the degree of $a$.  Then $n\geq 3$ and $G$ has $n$ vertices which are neighbors of $a$.  Since $G$ contains no triangles, none of these neighbors are adjacent to one another.   If $G$ has no vertices that are not neighbors of $a$ and $c$, then $G=K_{2,n}$.  Observe that the $2$-coloring of $K_{2,3}$ on the left in Figure~\ref{K23} has no non-trivial color preserving automorphisms.  Thus if $G=K_{2,n}$, then we must have $n\geq 4$.

\begin{figure}[htb]
\centerline{\includegraphics[width=4.25in]{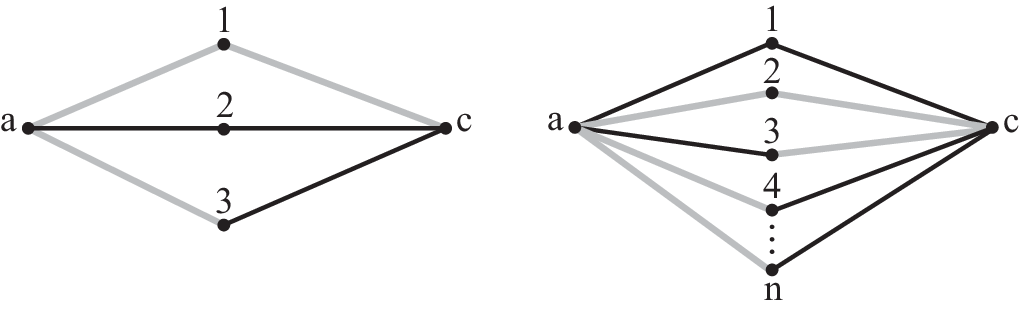}}
\caption{Intrinsically asymmetric 2-colorings of $K_{2,3}$ and $K_{2,n}$ for $n\geq 5$}
\label{K23}
\end{figure}

Now suppose that  $G=K_{2,n}$ and $n\geq 5$.  In this case, we color every edge containing vertex $a$ grey except for $\overline{a1}$ and $\overline{a3}$ which are colored black, and we color every edge containing vertex $c$ black except for $\overline{c2}$ and $\overline{c3}$ which are colored grey, as illustrated on the right in Figure~\ref{K23}.  Since $a$ is the only vertex of $K_{2,n}$ which is adjacent to at least three grey edges, every color preserving automorphism must fix $a$, and hence  also fix the edges $\overline{a1}$, $\overline{a2}$, and $\overline{a3}$.  Thus we can apply Lemma \ref{3edges} to conclude that this 2-coloring is intrinsically asymmetric in every orientable surface.  Hence if $G$ has no additional vertices that are not neighbors of $a$ and $c$, then $G=K_{2,4}$.

Next suppose for the sake of contradiction that $G$ contains a vertex $b$ which is not a neighbor of either $a$ or $c$.  Since $G$ is connected, without loss of generality $b$ is a neighbor of vertex $1$.  Now we can repeat the above argument with $b$ in place of $c$ to see that $b$ has the same set of neighbors as $a$ and $c$.  If $G$ has no additional vertices, then $G=K_{3,n}$ with $n\geq 3$, which we will see below is impossible.

First suppose that $G=K_{3,3}$ with partite sets $\{a,b,c\}$ and $\{1,2,3\}$.  We color $\overline{a3}$ and $\overline{b1c}$ black, and color the remaining edges of $G$ grey.  Now $h(1)=1$, since $1$ is the only vertex incident to two black edges.  Thus $h(a)=a$, $h(2)=2$, and $h(3)=3$. But now Lemma \ref{3edges} implies that this $2$-coloring of $G$ is intrinsically asymmetric in any orientable surface in which it embeds.  Thus we cannot have $G=K_{3,3}$.

Next suppose that $G$ contains $K_{3,4}$.  We label the partite sets of $K_{3,4}$ as $\{a,b,c\}$ and $\{1,2,3,4\}$.  Now we color $\overline{1a3c4}$ and $\overline{a2b}$ black and color all of the other edges in $G$ grey, and let $f$ be a color preserving automorphism of $G$.  Then $f$ fixes vertices $1$, $2$, $3$, and $a$.  But now Lemma \ref{3edges} implies that this $2$-coloring of $G$ is intrinsically asymmetric in any orientable surface in which it embeds.   Thus $G$ cannot contain $K_{3,4}$.  It now follows that $G\not =K_{3,n}$ for any $n\geq 3$.

 Therefore, if $G$ contains a vertex $b$ which is not a neighbor of $a$ or $c$, then it must also contain a vertex $d$ which is not a neighbor of $a$ or $c$.  Then by the above argument, $a$, $b$, $c$, and $d$ all have the same set of $n\geq 3$ neighbors.  But in this case, $G$ would contain a $K_{3,4}$ with partite sets $\{a,b,c,d\}$ and $\{1,2,3\}$. Since we saw above that this cannot occur, we must have $G=K_{2,4}$. \end{proof}\medskip
 
 The following lemma will be used in the proofs of both Theorem~\ref{SED(G)} and Theorem~\ref{exceptions}.

\begin{Lemma}\label{K1nSnm}  Suppose that $G$ is either $K_{1,n}$ with $n\geq 3$ or $S_{n,m}$ with $n$ and $m$ odd and at least one of $n$ or $m$ greater than $1$.  Then for any $2$-coloring of $G$ and any orientable surface $S$, there is an embedding $\Gamma$ of $G$ in $S$ and an orientation reversing homeomorphism of $(S,\Gamma)$  which induces a non-trivial color preserving automorphism of $\Gamma$.  \end{Lemma}

\begin{figure} [htb]
\centerline{\includegraphics[width=2.8in]{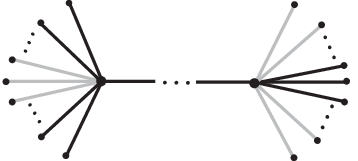}}
\caption{An orientation reversing color preserving homeomorphism interchanges the top and bottom edges.}
\label{star}
\end{figure}

\begin{proof}  First we consider $G=K_{1,n}$ with $n\geq 3$, and suppose that $G$ is $2$-colored and $S$ is an orientable surface.  Now let $\Gamma$ be an embedding of $G$ in $S$ such that the grey edges are grouped together on the left and the black edges are grouped together on the right as illustrated in Figure \ref{Sn}.  Then there is an orientation reversing color preserving homeomorphism $h$ of $(S,\Gamma)$ which interchanges the grey edges on the top with those on the bottom and interchanges the black edges on the top with those on the bottom.  If the number of edges of a given color is odd then $h$ will fix one such edge.  If there is an even number of edges of a given color, then $h$ fixes none of the edges of that color.  
 
\begin{figure} [h]
\centerline{\includegraphics[width=1.2in]{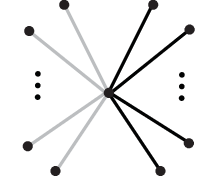}}
\caption{An orientation reversing color preserving homeomorphism interchanges the top and bottom edges.}
\label{Sn}
\end{figure}

Next we consider $G=S_{n,m}$ with both $n$ and $m$ odd and at least one of $m$ or $n$ greater than $1$, and suppose that $G$ is $2$-colored and $S$ is an orientable surface.  We will refer to edges with a vertex of degree 1 as ``pendant'' edges. Then the number of pendant edges incident to each vertex is odd.  Thus at each vertex, there must be an odd number of pendant edges of one color and an even number (possibly zero) of pendant edges of the other color.  Now let $\Gamma$ be an embedding of $S_{n,m}$ in $S$ such that at a given vertex the odd number of pendant edges of a single color are in the center and the remaining pendant edges are divided evenly between the top and the bottom as illustrated in Figure~\ref{star}.  Then there is an orientation reversing color preserving homeomorphism of $(S, \Gamma)$ which interchanges the top and bottom pendant edges of $\Gamma$.  \end{proof}

\medskip

We are now ready to prove Theorem \ref{SED(G)}.

\begin{Theorem} A graph $G$ has $SED(G) > 2$ if and only if $G$ is a single edge, $C_3$, $C_4$, $C_5$, $K_4$, $K_5$, $K_{2,4}$, $K_{1,n}$ with $n\geq 3$, or $S_{n,m}$ with $n$ and $m$ odd and at least one of $n$ or $m$ greater than $1$.  \end{Theorem}

\begin{proof}  We begin by showing that none of the graphs listed in the theorem has $SED(G)\leq 2$.  Certainly, if $G$ is a single edge, then $SED(G)\not\leq 2$.  Every non-trivial $2$-coloring of $C_3$, $C_4$, and $C_5$ (up to switching black and grey) is illustrated in Figure~\ref{CnChart}.  It is easy to check that for any embedding $\Gamma$ of each of these $2$-colored graphs in an orientable surface $S$, there is a homeomorphism of $(S,\Gamma)$ inducing a non-trivial color preserving automorphism on $\Gamma$.  Thus if $G$ is any of these graphs then $SED(G)>2$.

\begin{figure}[h!]
\centerline{\includegraphics[width=.70\textwidth]{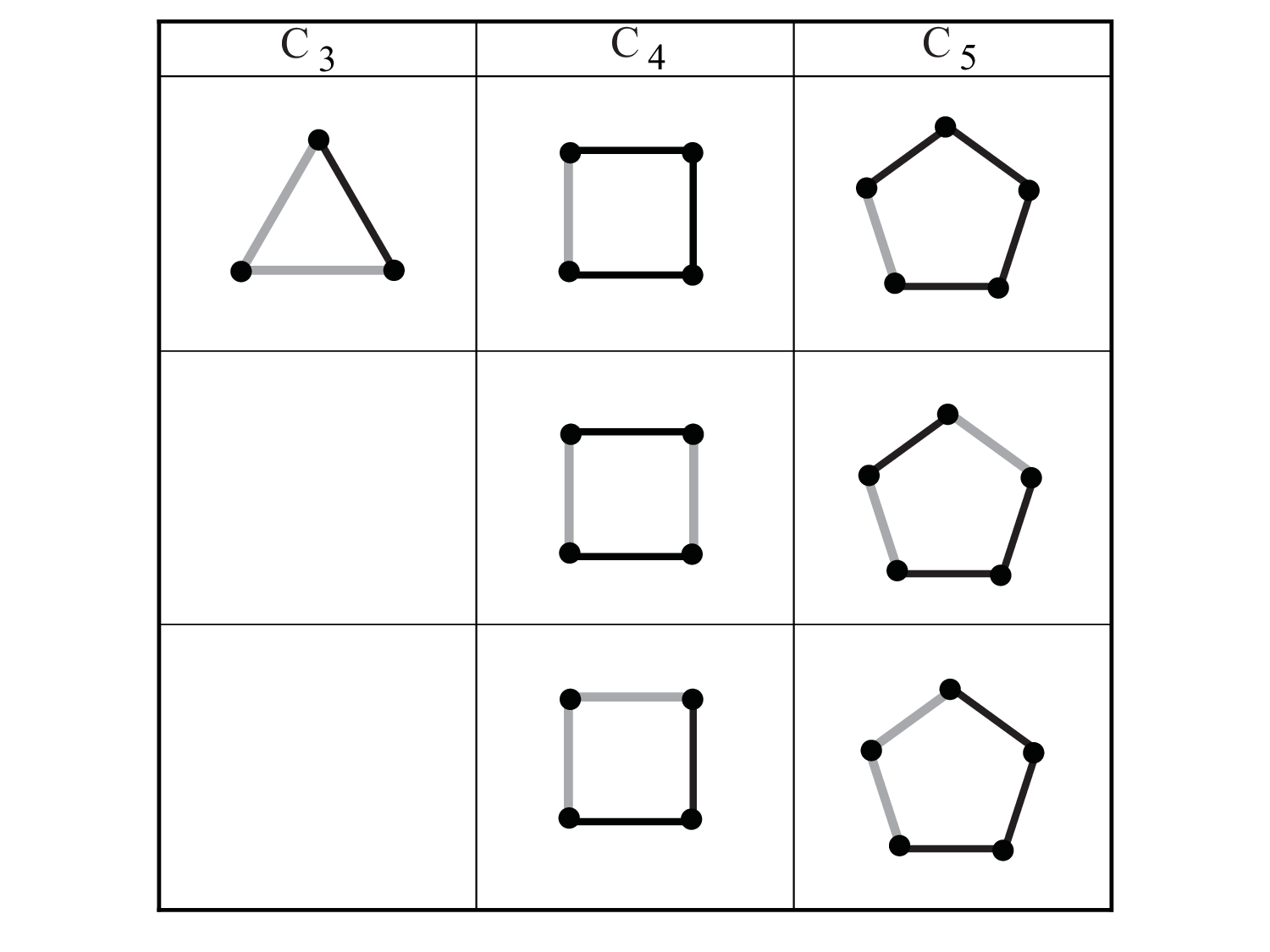}}
\caption{Every $2$-coloring of $C_3$, $C_4$, and $C_5$.}
\label{CnChart}
\end{figure}

For the graph $K_4$, the chart in Figure~\ref{K4Chart} illustrates all possible $2$-colorings.  Furthermore, in the proof of Theorem~\ref{K4Theorem} we showed that each of the $2$-colorings in the first two rows of the chart is induced by a reflection of an embedding of $K_4$ in $S^2$, and the $2$-coloring in the third row of the chart is induced by a reflection of the embedding of $K_4$ in a torus illustrated in Figure~\ref{K4T2}.  Thus $SED(K_4)>2$.

Now we consider the graph $K_5$.   In Figure~\ref{K5Chart} we illustrate all subgraphs (including those which are disconnected) of $K_5$ with up to $5$ edges.  Without loss of generality, for any non-trivial $2$-coloring of $K_5$ the black (possibly disconnected) subgraph is one of those illustrated.

\begin{figure}[h!]
\centerline{\includegraphics[width=2.6in]{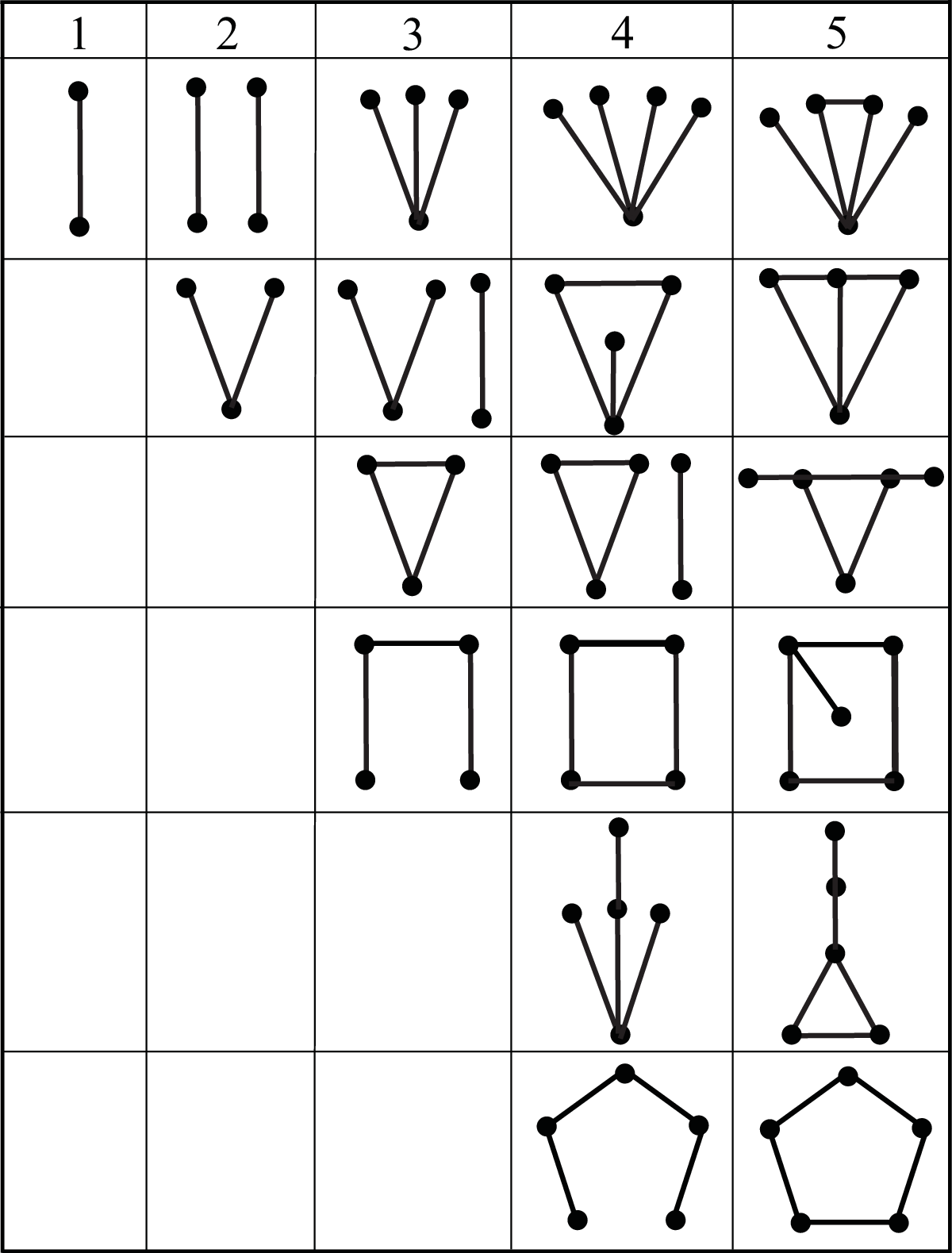}}
\caption{All subgraphs of $K_5$ with up to $5$ edges.}
\label{K5Chart}
\end{figure}

Each of these $2$-colorings of $K_5$ has at least one automorphism of order $2$.  For each coloring, we choose a labeling of the vertices of $K_5$ such that the automorphism $(24)$ or $(23)(14)$ is color preserving.  Now we embed $K_5$ in a torus as illustrated in Figure~\ref{K5Torus}, and observe that the automorphisms $(24)$ and $(13)(24)$ are each induced by a reflection of the torus. It follows that $SED(K_5) > 2$.

\begin{figure}[h!]
\centerline{\includegraphics[width=1.5in]{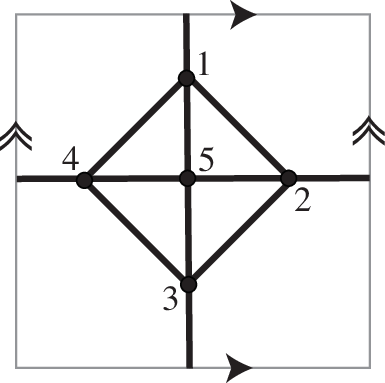}}
\caption{$(24)$ and $(23)(14)$ are induced by reflections of this embedding of $K_5$ in a torus.}
\label{K5Torus}
\end{figure}

 Next we consider $2$-colorings of $K_{2,4}$.  Let $v$ and $w$ be the vertices of valence $4$ and let $1$, $2$, $3$, and $4$ be the vertices of valence $2$.  Suppose that for some $i\not= j$, the oriented paths $viw$ and $vjw$ have identical colorings.  Hence $K_{2,4}$ is colored as illustrated on the left in Figure~\ref{K24Reflection}.  Then we embed $K_{2,4}$ in a sphere as illustrated  on the right in Figure~\ref{K24Reflection}, and observe that a reflection through the equator induces the automorphism $(12)$.

\begin{figure} [htb]
\centerline{\includegraphics[width=4in]{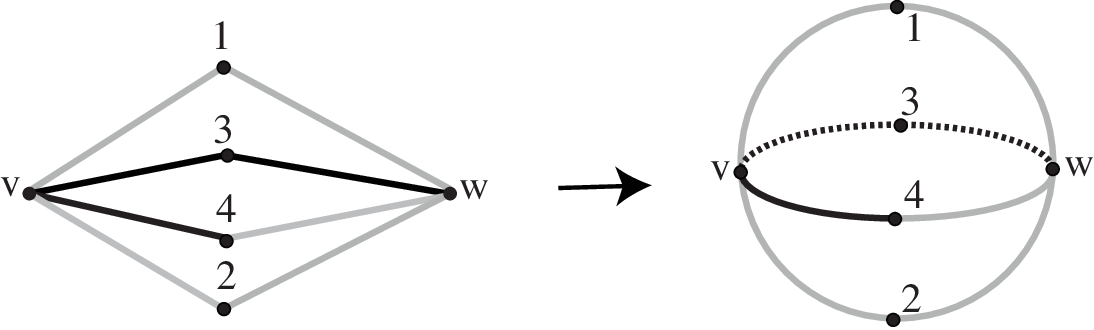}}
\caption{A reflection through the equator induces the automorphism $(12)$.}
\label{K24Reflection}
\end{figure}

Now suppose that no pair of paths from $v$ to $w$ have the same coloring.  Hence $K_{2,4}$ is colored as illustrated on the left in Figure~\ref{K24}.  We embed $K_{2,4}$ in a torus as illustrated on the right in Figure~\ref{K24}.  Then the  homeomorphism of the torus obtained by composing a reflection through the circle $1v2w$ with a translation along the circle induces the automorphism $(vw)(12)$.  Hence $SED(K_{2,4})>2$.

\begin{figure} [h!]
\centerline{\includegraphics[width=4in]{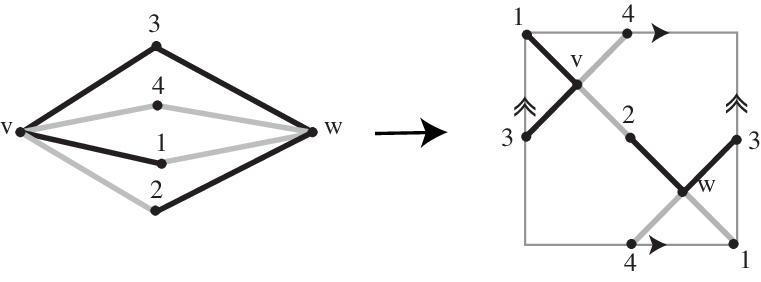}}
\caption{The automorphism $(vw)(12)$ is induced by an orientation reversing homeomorphism of a torus.}
\label{K24}
\end{figure}

Finally, by Lemma~\ref{K1nSnm}, if $G$ is either $K_{1,n}$ with $n\geq 3$ or $S_{n,m}$ with $n$ and $m$ odd and at least one of $n$ or $m$ greater than $1$, then $SED(G)>2$.  Thus for each of the graphs $G$ listed in the theorem, we have shown that $SED(G)>2$.

\medskip

In order to prove the converse, let $G$ be a graph with $SED(G)>2$.  First suppose that $G$ does not have a vertex of degree at least $3$.  Then $G$ must be a cycle $C_n$ or a path.  If $G$ is a path whose length is greater than one, then we can color an edge at one end black and the rest of the edges of $G$ grey to get a $2$-coloring which is intrinsically asymmetric in every orientable surface.  Thus if $G$ is a path, then $G$ must be a single edge.

Suppose that $G=C_n$ and $n\geq 6$.   Label the vertices of $C_n$ consecutively as $1$, $2$, \dots, $n$.  Now color $\overline{123}$ and $\overline{45}$ black and color the rest of the edges grey.  Since vertices $3$ and $4$ are adjacent but $5$ and $1$ are not, any color preserving automorphism of $G$ fixes vertices $1$, $2$, $3$, $4$, and $5$.  Thus $G$ has no non-trivial color preserving automorphism.  Hence again $ED(G) = 2$ contradicting our hypothesis. It follows that if $G$ is a cycle then $G$ is $C_3$, $C_4$, or $C_5$.

Next suppose that $G$ has a vertex of degree at least $3$.  If $G$ contains a triangle, then by Lemma \ref{RefereeLemma1}, $G$ is either $K_4$ or $K_5$.  So we assume that $G$ does not contain a triangle.  If some vertex of degree at least $3$ is adjacent to at least two vertices whose degrees are at least $2$, then it follows from Lemma \ref{RefereeLemma2} that $G=K_{2,4}$.  Thus we assume that every vertex of degree at least $3$ has at most one neighbor whose degree is at least $2$.  

Now suppose that $G$ has precisely one vertex $v$ of degree at least $3$. If all of the neighbors of $v$ have degree $1$, then $G=K_{1,n}$ for some $n\geq 3$.  If precisely one neighbor of $v$ has degree $2$, then $G$ is $S_{1,n}$ for some $n>1$.   Thus we assume that $G$ has at least two vertices $u$ and $v$ of degree at least $3$. Since $G$ is connected, $G$ contains a simple path $P$ joining $u$ and $v$.  Since $u$ and $v$ each have at most one neighbor whose degree is at least $2$, all of the neighbors of $u$ and $v$ not contained in $P$ have degree $1$. If $G$ contains a vertex distinct from $u$ and $v$ whose degree is at least $3$, then $P$ contains a vertex $w$ of degree at least $3$ in its interior.  It now follows that $w$ is adjacent to a pair of vertices whose degrees are at least $2$.  As this is contrary to the assumption at the end of the above paragraph, $u$ and $v$ are the only vertices of degree at least $3$.  Thus $G=S_{n,m}$, where $n>1$ and $m>1$.  

It remains to show that if either $G=S_{1,n}$ or $G=S_{n,m}$ where $n>1$ and $m>1$, then $n$ and $m$ must be odd.  Suppose for the sake of contradiction that at least one of $n$ or $m$ is even.  Without loss of generality, $n\geq 2$ is even and the vertex $v$ has $n$ neighbors of degree $1$ and one neighbor of degree more than $1$.     Let $e_1$ be a pendant edge incident to $v$.  We color $e_1$ black and color all of the remaining edges of $G$ grey.  Now suppose $\Gamma$ is an embedding of $G$ in an orientable surface $S$, and let $h$ be a color preserving homeomorphism of $(S, \Gamma)$. Then $h(e_1)=e_1$ and $h(v)=v$.  Since there are an odd number of grey pendant edges incident to $v$, at least one such edge $e_2$ must be fixed by $h$.  Finally, the only edge incident to $v$ which is not pendant must also be fixed by $h$.  Thus by Lemma~\ref{3edges}, this $2$-coloring of $G$ is intrinsically asymmetric in $S$.  But since $S$ was arbitrary this $2$-coloring is intrinsically asymmetric in every orientable surface, contradicting our hypothesis that $SED(G)>2$.  Hence if $G=S_{n,m}$, then both $n$ and $m$ must be odd.
\end{proof}
\medskip


\section{Graphs with Intrinsically Chiral $2$-Colorings}


We now prove Theorem 3 which characterizes those graphs with a $2$-coloring that is intrinsically chiral in every orientable surface.  Note that if a graph $G$ does not embed in an orientable surface $S$, then $G$ is vacuously intrinsically chiral in $S$.

\setcounter{Theorem}{2}
\begin{Theorem}  A graph $G$ has a $2$-coloring which is intrinsically chiral in every orientable surface if and only if $G$ has at least one vertex of degree at least $3$ and $SED(G)\leq 2$.  \end{Theorem}

\begin{proof}  First suppose that $G$ has no vertex of degree at least $3$.  Then $G$ can be embedded in a circle.  Hence regardless of how the edges of $G$ are colored, $G$ can be embedded in any orientable surface so that it is pointwise fixed by a reflection of the surface.  Thus no $2$-coloring of $G$ can be intrinsically chiral in any orientable surface.  

Thus we suppose that $G$ has at least one vertex of degree at least $3$ and $SED(G)>2$.  The proof of Theorem~\ref{SED(G)} shows that for every $2$-coloring of $K_4$, $K_5$, and $K_{2,4}$, there is an embedding  $\Gamma$ of the colored graph in a sphere or a torus such that some orientation reversing homeomorphism of the surface induces a color preserving automorphism on the graph.  Thus no $2$-coloring of these graphs can be intrinsically chiral in every orientable surface.

 Now since $SED(G)>2$, by Theorem~\ref{SED(G)}, $G$ must be either $K_{1,n}$ with $n\geq 3$ or $S_{n,m}$ with both $n$ and $m$ odd and at least one of $n$ or $m$ greater than $1$.  However, by Lemma~\ref{K1nSnm}, for any $2$-coloring of $G$ and any orientable surface $S$, there is an embedding $\Gamma$ of $G$ in $S$ and an orientation reversing homeomorphism of $(S,\Gamma)$  which induces a non-trivial color preserving automorphism of $\Gamma$.  Thus again no $2$-coloring of $G$ can be intrinsically chiral in any orientable surface.  

\medskip

To prove the converse, suppose that $G$ has at least one vertex whose degree is at least $3$ and $SED(G)\leq 2$.   Since $SED(G)\leq 2$, there is a $2$-coloring of $G$ which is intrinsically asymmetric in every orientable surface.  Now suppose for the sake of contradiction that $\Gamma$ is an embedding of this $2$-coloring of $G$ in an orientable surface $S$ such that there is a color preserving orientation reversing homomorphism $h$ of $(S,\Gamma)$.  Since this $2$-coloring is intrinsically asymmetric in $S$, $h$ induces the trivial automorphism on $\Gamma$.  

Now there is a disk neighborhood $D$ of a vertex $v$ of $\Gamma$ such that $h(D)=D$ and $D$ contains no vertices other than $v$.  Since the boundary of $D$ intersects at least three edges, each of which is fixed by $h$, we know that $h$ preserves the orientation of $D$.  But this is contrary to our assumption that $h$ was orientation reversing on $S$.  Hence $G$ must be intrinsically chiral in all orientable surfaces. \end{proof}

\medskip

\section{Acknowledgements} The authors are grateful to an anonymous referee whose suggestions enabled us to strengthen our results while simplifying our proofs.

\end{document}